# Robust Vehicle Routing Problem in the Presence of a Proactive Attacker

H. Bigdeli[1], S.M.S. Mirdamadi[2], Adrian Deaconu[3] , J. Tayyebi[4]

1. Department of Science and Technology Studies, AJA Command and Staff University, Tehran, Iran.

2. Department of Industrial Engineering, Isfahan University of Technology, Isfahan, Iran

3. Department of Mathematics and Computer Science, Transilvania University, Brasov, 500091, Romania.

4. Department of Industrial Engineering, Birjand University of Technology, Birjand, Iran.

**Abstract:**

The vehicle routing problem has great importance and application in transportation and supply chain management. In this case, there are several supply requests in a transportation network. The main goal is to allocate customers to available vehicles and find the sequence of customer visits on each route. It is possible to attack the arcs of the network when two actors compete against each other with conflicting goals. In this case, the attacker may destroy the network arcs or prevent the supply of customers. The defender also intends to fulfill customer demand with the least cost and the most success rate. To this end, we define the robust vehicle routing problem in the presence of a proactive attacker. In this research, we propose a mathematical model with two objective functions to minimize total cost and maximize the total success rate in fulfilling the demand. Also, an exact algorithm for generating feasible routes has been used to solve small-scale problems. To solve large-scale problems, a simulated annealing algorithm is used. We test one numerical example problem and six medium and large size problems to evaluate the proposed algorithms' efficiency. The results demonstrate that the total success rate in supply is directly related to the amount of cost; crossing shorter routes decrease the total of success rate.

## 1. Introduction

The vehicle routing problem (VRP) is one of the most practical problems in supply chain management. It seeks to select and allocate possible routes to vehicles from distribution centers to the customers to minimize the associated costs. In this regard, the optimal allocation of vehicles to customers and routing have a significant effect on reducing costs. However, satisfying all demands is difficult due to critical situations such as relief logistics or military environments. In this situation, due to the insecurity of the routes or the possibility of route attack, it will be probable not to cover all demands.

Vehicle Routing is a practical problem in the field of combinatorial optimization. There is a lot of research in the field of logistics and supply chain associated with the vehicle routing problem. The goal is to distribute enough goods between customers and depots through a set of routes to satisfy all customers, when an objective function (such as total traveled distance or total cost) is optimized. In this case, all demands must be satisfied and the total demand for each route should not exceed the capacity of the vehicle. Solving VRP involves allocating customers to routes and then determining the sequence of customer visit. The basic VRP problem is the capacitated vehicle routing (CVRP), in which a set of homogeneous fleets deliver goods from the depot to the customers.

In the VRP, the network of roads is considered as a graph. The arcs represent the connections and the nodes indicate the location of customers and the depot. The traveling salesman (TSP) problem is the basis of the vehicle routing problem. It seeks to find a route that passes through all nodes and minimizes the total cost of passing. The multiple traveling salesman problem (MTSP) is an extension of the traveling salesman problem with unlimited capacity. If limited capacity is considered for each of the vehicles, it converted into VRP. The extensive literature on VRP and its variants can be found in Cordeau et al. [1], Laporte [2],Golden et al. [3], and Toth and Vigo [4].

In the ground vehicle routing problem, there is a set of transport requests and a set of transport fleets in a network. The goal is to fulfill the demand of network nodes with objectives such as minimizing cost or maximizing demand coverage. Also, due to the unstable situation, there is a possibility of attacking network. Therefore routes must be selected to achieve the goals with the existence of unforeseen attacks.

In VRP, there are a network and several vehicles with a certain capacity at a specific point, and the goal is to find proper routes according to the objective function. Figure 1 depicts an example of a vehicle routing problem.

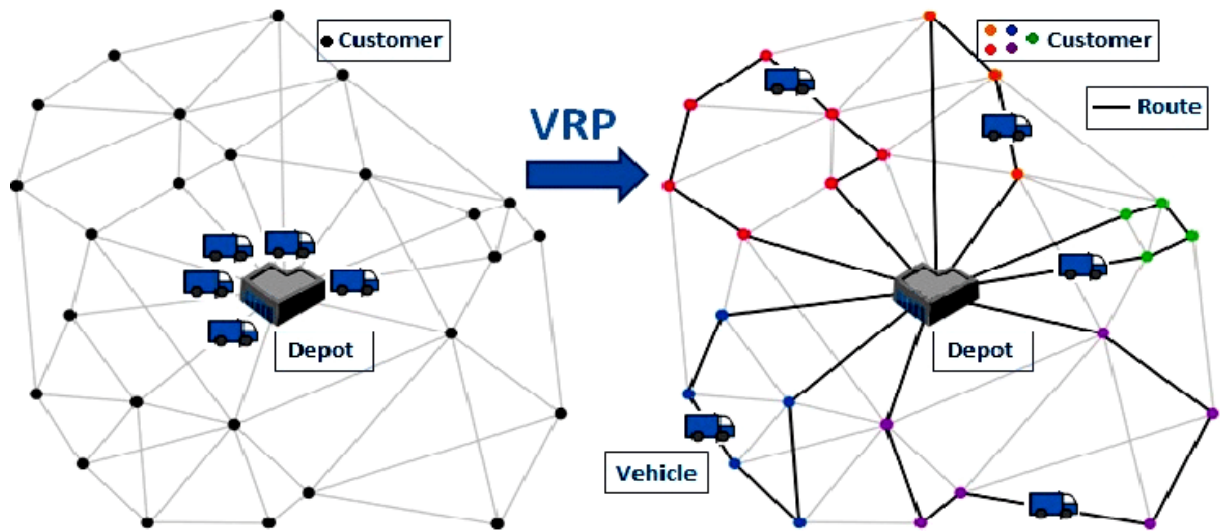

Figure 1. A VRP example

## 2. The Classical Capacitated Vehicle Routing Problem

In the Capacitated Vehicle Routing Problem, a point is known as the central warehouse or depot and is usually displayed as a zero point. There are also n demand points, which are considered as a set of demand points $N = \{1,2,...,n\}$. The amount of demand for each point $i \in N$ is the amount of goods that must reach the desired point and is denoted by $d_i$. There is also a fleet which is a set of homogeneous transport vehicles in the central depot $K = \{1,2,...,|K|\}$. Each of these homogeneous vehicles has a definite capacity of $Q > 0$. Each vehicles serves a subset of points $S \subseteq N$ that start from the depot point and then visit all points. $c_{ij}$ shows the cost of traveling on the arc $(i,j)$.

As mentioned, the problem can be considered as a graph. $V = \{0,1,2, ... n\}$ determines the nodes and $E = \{(i,j) \mid i, j \in V, i \neq j\}$ is the set of arcs or transport paths. Therefore, the graph $G = (V, E)$ indicates the road network.

In VRP, a route is a sequence of points that are represented by $r = (i_0, i_1, ..., i_s, i_{s+1})$. The beginning and end point of the route are the same depot ($i_0 = i_{s+1} = 0$). The route cost is also calculated as $c(r) = \sum_{p=0}^{s} c_{i_p, i_{p+1}}$. A route is called feasible if the total demand that each vehicle carries is less or equal than the available capacity, i.e., $\sum_{i \in S} d_i \leq Q$. Also, each customer should be visited only once. Therefore, to solve the CVRP, two problems must be resolved:

- Independent categorization of $N$ customers in $|K|$ separate sets so that capacity limit is satisfied.
- Routing each vehicle from the depot to the customer points and return to the depot.

The second part of problem i.e., routing between nodes and returning to the depot, is the same as the traveling salesman problem [5]. It should be noted that both parts of the vehicle routing problem are intertwined. In fact, the answer of each part affects the other part.

According to the descriptions, the mathematical model of the classical capacitated vehicle routing problem is the following formulations 1-7.

$$Min\ Z = \sum_{(i,j) \in E} \sum_{k \in K} c_{ij} x_{ijk} \tag{1}$$

s.t.

$$\sum_{i \in V} \sum_{k \in K} x_{ijk} = 1 \qquad \forall j \in V \backslash \{0\} \tag{2}$$

$$\sum_{j \in V} x_{jik} = \sum_{j \in V} x_{ijk} \qquad \forall i \in V, \forall k \in K \tag{3}$$

$$u_{ik} = 0 \qquad \forall k \in K, i = 0 \tag{4}$$

$$u_{ik} + 1 \leq u_{jk} + M.(1 - x_{ijk}) \qquad \forall (i,j) \in E, j \neq 0, \forall k \in K \tag{5}$$

$$\sum_{(i,j) \in E} d_i x_{ijk} \leq Q_k \qquad \forall k \in K \tag{6}$$

$$x_{ijk} \in \{0,1\}, u_{ik} \geq 0 \qquad \forall (i,j) \in E, \forall k \in K \tag{7}$$

Equation (1) represents the cost or distance traveled by the homogeneous fleet to fulfill demand at different points. Each customer must be serviced by a vehicle that expresses in constraint (2). Constraint (3) expresses the flow balance in the network. In other words, if a vehicle enters a node, it must leave that node. Equations (4) and (5) are used to sub-tour elimination in the problem. Therefore, all created routes should start from the depot and return to the depot after visiting a few customers. Constraint (6) is defined to comply with vehicle capacity. Finally, Equation (7) expresses the type of the variables.

## 2.1. Numerical example

Consider a problem with eight demand points and a central warehouse. In this example, there are two vehicles with a capacity of 180 in the depot. Demand values and point coordinates are represented in Table 1.

Table 1. Sample data for the capacitated vehicle routing problem

| Node number | x coordinate | y coordinate | Demand ($d_i$) |
|---|---|---|---|
| depot(0) | 7 | 47 | 0 |
| 1 | 28 | 45 | 47 |
| 2 | 21 | 1 | 35 |
| 3 | 38 | 10 | 47 |
| 4 | 48 | 0 | 33 |
| 5 | 44 | 31 | 28 |
| 6 | 17 | 35 | 43 |
| 7 | 12 | 47 | 39 |
| 8 | 3 | 12 | 36 |

The mathematical model of the capacitated vehicle routing problem is implemented in GAMS 24.1 software using the CPLEX solver. The cost function is calculated based on the distance of the nodes according to their coordinates. After solving the problem, the optimal solution consists of two separate routes. A summary of the results is provided in Table 2. An overview of the optimal route depicted in Figure 2.

Table 2. Details of the solution to the sample CVRP

| route number | route | used capacity | available capacity | objective function | time (s) |
|---|---|---|---|---|---|
| 1 | 0→8→2→4→3→5→0 | 179 | 180 | 211.25 | 3.86 |
| 2 | 0→6→1→7→0 | 129 | 180 | | |

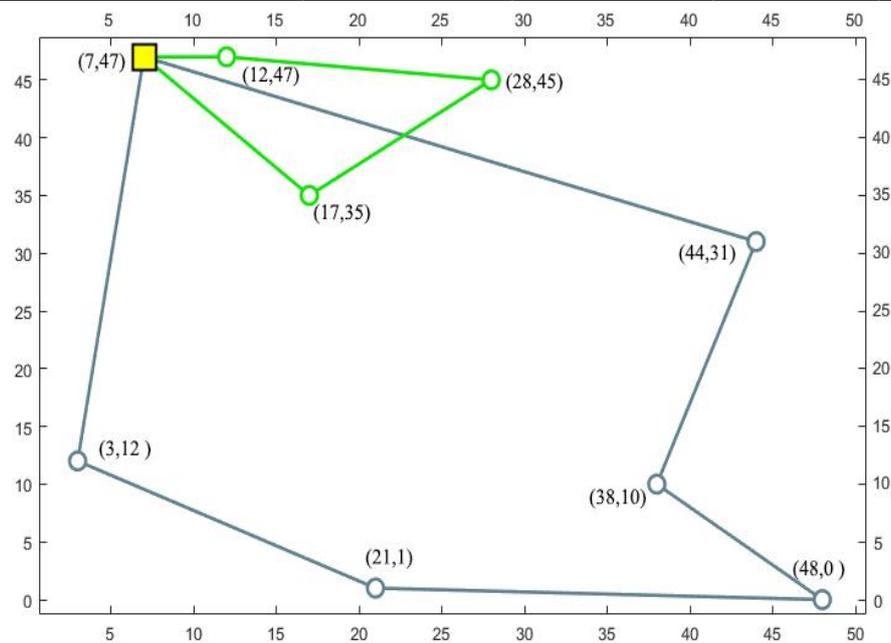

Figure 2. Representation of the optimal solution

## 2.2. Route-based mathematical modeling

The extensive CVRP formulation was proposed by Balinski and Quandt [6]. Their model was presented based on set partitioning model. In this case, feasible routes are generated, then the mathematical model selects the generated routes. Thus, the mathematical model selects the best routes from the generated feasible routes.

Assume that $r \in \Omega$ represents all feasible routes of CVRP. Here, $r = (i_0, i_1, \ldots, i_s, i_{s+1})$ where $i_0$ and $i_{s+1}$ are depot point. The cost of each route is calculated as $c_r = \sum_{j=0}^{s} c_{i_j, i_{j+1}}$. The coefficient $a_{ir} \in \{0,1\}$ is equal to one if customer $i$ is visited in the route $r$. On the other hand, the variables $\lambda_r$ is equal to one if $r$ is selected, otherwise, it is zero. The total number of fleet vehicles in the problem is also equal to $|K|$. The problem is formulated as follows:

$$Min\ Z = \sum_{r \in \Omega} c_r \lambda_r \tag{8}$$

s.t.

$$\sum_{r \in \Omega} a_{ir} \lambda_r = 1 \qquad \forall i \in N \tag{9}$$

$$\sum_{r \in \Omega} \lambda_r \leq K \tag{10}$$

$$\lambda_r \in \{0,1\} \tag{11}$$

Equation (8) shows the objective function and represents the cost of the selected routes in the optimal solution. Each customer should be visited only once, as stated in equation (9). Constraint (10) determines the number of available fleets. Equation (11) determines the type of variables.

## 3. Robust Vehicle Routing Problem in the Presence of a Proactive Attacker

The vehicle routing problem was first addressed by Dantzing and Ramser [7]. After that, VRP became one of the most important and practical problems in logistics and combinatorial optimization. The main purpose is to fulfill demand and minimize operational costs, which are usually equal to the distance traversed by all vehicles. There are various constraints in the vehicle routing problem, e.g., capacity constraint and the number of available fleets.

In the face of natural or man-made disasters, proper distribution and supply of vital goods are very important. These include distribution during events such as war or natural disasters such as earthquakes. Therefore, many efforts to save human lives are based on reducing human suffering, reducing the cost of providing essential commodities, and reducing the risk of loss of human life. In normal circumstances, routing is used to satisfy the required demand, but when a vehicle is likely to break down or not reach the destination, it will be more difficult to routing well. In natural disasters, this can be caused by aftershocks or dangerous routes. While in wartime, destructive enemy actions or roadside mines destroy roads or vehicles.

In this research, the robust vehicle routing problem in the presence of a proactive attacker is discussed. The problem can be used in humanitarian logistics and military applications. In this case, it is assumed that each arc has a certain chance of being prevented and destroyed. There are two attacking forces with conflict goals. Defender intends to fulfill the demand for several points in the network. On the other hand, the attacker tries to prevent the actions of the defender by destroying or attacking the arcs of the network. The probability of the attacking arcs can be estimated by the defender. When an arc is attacked, the vehicle can't move through it, so the vehicle will not be able to continue the road or can't fulfill the demand. Also in business operations, the supply fleet may be damaged during transportation or prevented from reaching the destination due to route traffic. Therefore, the routes used in most routing problems are usually unstable and there is a possibility of attacking. The purpose of evaluating this problem is to investigate the impact of different risks on the VRP. Based on this, the robust vehicle routing problem in the presence of a proactive attacker is proposed.

In this study, the possibility of fulfillment demand is not considered definitively. Factors such as the presence of offensive forces or natural disasters result in uncertainty in fulfillment or failure to reach the destination. In fact, an estimated amount of "attacking probability" is assumed for the arcs. The value of the probability can be predicted or estimated by experts or based on various environmental data. Partial destruction, complete destruction, or blockage of roads are types of roadblocks. When a route is blocked, the vehicle will not be able to fulfill the demand.

Two different modes are considered to fulfill customer demand. In the first case, the amount of inventory is enough to fulfill all the demand and the goal is to minimize costs. In the latter case, the amount of inventory is less than the required amount and the goal is to increase the probability of demand fulfillment given the risks involved. This situation usually occurs during the midst of disaster or wartime, because the available resources are reduced and the demands are increased. Therefore, in this research, two issues are examined. The first is related to reducing costs and the other is related to maximizing fulfilled demand.

Various types of uncertainty have been studied in the literature on vehicle routing problems, but uncertainties in the route or fulfillment have not been addressed. Three types of uncertainty that are widely considered in the literature are stochastic demand, stochastic customers, and stochastic travel time. Various studies on the vehicle routing problem have been reviewed and categorized by Berhan et al. [8].

- The VRP with stochastic demand was analyzed in 1969 [9]. In this case, the demand of each node is uncertain. Meta-heuristic approaches are among the applied research in this field [10].
- The second category is related to VRP with uncertain customers. Therefore, the presence or absence of the customer in a node, the weather during distribution, or the possibility of distribution are considered uncertain [11].
- The third category includes the VRP with stochastic travel service time. In this class, service time and travel time are considered uncertain [12].

In the stochastic VRP field, there are very little research has reported vehicle loss, destruction, or inability to fulfill the demands. For example, the vehicle routing problem, in which there was a possibility of failure of the transport fleet, was investigated by Liu et al. [13]. But attacking does not cause destruction or vehicle loss. Instead, they are rerouted based on a two-stage model.

The concept of probability of attacking is adapted from the interdiction network problem. So, the conflict between the leader and the follower is considered. In this model, goals such as finding the shortest path or maximum flow are considered. However, the attacker has pursued interdiction of arcs or nodes to block and disrupt the defender's movements [14].

## 4. Mathematical Formulation

Suppose a graph $G = (V, E)$ for the vehicle routing problem with the possibility of attacking of arcs. $V$ represents the set of nodes and $E$ represents the set of arcs in the network. Each node has a definite demand $d_i$. The cost of passing the arc $(i, j)$ is equal to $c_{ij}$. There may also be an attacking probability, which will prevent the demand fulfillment or will damage the vehicle. The attacking probability of arc $(i, j)$ is defined based on a $poi_{ij}$ parameter ($0 \leq poi_{ij} \leq 1$). Therefore, the probability of demand fulfillment success is equal to $p_{ij} = 1 - poi_{ij}$. The probability of attacking or success in crossing an arc is estimated values that can be predicted according to available data. Also, these values can be suggested by experts according to environmental conditions and data.

The mathematical model presented in this section is a route-based formulation. Assume that all possible routes are represented as a set $\Omega$. Each route $w \in \Omega$ is defined based on the nodes that are visited by vehicle $w$. The cost of each route can be defined as $c_w = \sum_{(i,j) \in w} c_{ij}$. If any route $w$ is defined as $w = \{0, n_1, \ldots, n_w, 0\}$. Then the probability of success in demand fulfillment for $n_1$ is equal to $\square_{n_1 w} = p_{0n_1}$. Also, the probability of demand fulfillment success for node $n_2$ is equal to $\square_{n_2 w} = p_{0n_1} p_{n_1 n_2}$. Thus, the parameter $\square_{iw}$ indicates the probability of reaching node $i$ with demand fulfillment success.

In this research, two different conditions are considered for demand fulfillment. In the first case, the available amount of goods is greater than the demand, and the main goal is travel cost minimization. In the second case, the available amount of goods is less than the demand, so the objective function is to maximize the probability of demand fulfillment.

In situations where the amount of supply exceeds the demand, the main goal is usually to reduce the associated costs and increase the efficiency of supply. According to the stated cases, the mathematical model is defined based on the sets, parameters, and variables of Tables 3-5. Due to the lack of knowledge of the optimal routes before finding the solution, it will be impossible to calculate the success rate $\square_{ik}$. The classical CVRP formulation is defined based on the arc, while to calculate the success rate, we must know the sequence of customers' visits. Therefore, the route-based formulation is used.

Table 3. Mathematical model sets

| Symbol | Description |
|---|---|
| $V$ | Set of nodes and i is the corresponding index $i \in V = \{0,1,2, \ldots, N\}$ |
| $E$ | Set of arcs $E = \{(i,j) \mid i, j \in V, i \neq j\}$ |
| $K$ | Set of available vehicles $K = \{1,2, \ldots, |K|\}$ |
| $\Omega$ | The set of all feasible routes is based on the constraints and $w$ is the corresponding index. |

Table 4. Mathematical model parameters

| Symbol | Description |
|---|---|
| $d_i$ | The amount of demand |
| $c_w$ | Amount of $wth$ route cost |

| | | |
|---|---|---|
| $\Phi_{iw}$ | The probability of success for demand fulfillment (probability of reaching point $i$ demand fulfillment success) | |
| $a_{iw}$ | if node $i$ is in route $w$; otherwise $a_{iw} = 0$ $a_{iw} = 1$ | |
| $Q_w$ | vehicle capacity | |
| $\alpha$ | Minimum probability of success for demand fulfillment | |
| $\beta$ | Maximum acceptable amount for total costs | |

Table 5. Mathematical model variables

| Symbol | Description |
|---|---|
| $x_w$ | if route $w$ is selected; otherwise $x_w = 0$ $x_w = 1$ |

The robust vehicle routing problem in the presence of a proactive attacker is described as 12-17 formulation in cost minimization form.

$$Min \ Z = \sum_{w \in \Omega} c_w x_w \tag{12}$$

s.t.

$$\sum_{w \in \Omega} a_{iw} x_w = 1 \qquad \forall i \in V/\{0\} \tag{13}$$

$$\sum_{k \in \Omega} a_{iw} x_w d_i \leq Q_w \qquad \forall w \in \Omega \tag{14}$$

$$\sum_{w \in \Omega} a_{iw} x_w \Phi_{iw} \geq \alpha \qquad \forall i \in V/\{0\} \tag{15}$$

$$\sum_{w \in \Omega} x_w \leq |K| \tag{16}$$

$$a_{iw}, x_w \in \{0,1\} \qquad \forall w \in \Omega \tag{17}$$

Equations (12-17) represent the cost minimization model. Equation (12) shows the objective function and seeks to minimize the total cost. Equation (13) fulfills all the demand node among the selected routes. The capacity limit of vehicles is defined by Equation (14). Each feasible route must have at least an $\alpha$ probability of success in demand fulfillment for route points, as constrained by (15). Hence, the probability of success in demand fulfillment should not be less than $\alpha$. Equation (16) shows the maximum number of vehicles. Equation (17) describes the type of variables.

When available goods are less than the total demand, it is better to increase the probability of success in demand fulfillment as much as possible. Therefore, the main goal is to maximize the demand fulfillment so that the attacker will not block the routes. As a result, the maximization formulation can be state as follow:

$$Max \ Z = \sum_{w \in \Omega} a_{iw} x_w \Phi_{iw} \tag{18}$$

s.t.

$$\sum_{w \in \Omega} a_{iw} x_w = 1 \quad \forall i \in V/\{0\} \quad (19)$$

$$\sum_{w \in \Omega} a_{iw} x_w d_i \leq Q_w \quad \forall w \in \Omega \quad (20)$$

$$\sum_{w \in \Omega} c_w x_w \leq \beta \quad \forall i \in V/\{0\} \quad (21)$$

$$\sum_{w \in \Omega} x_w \leq |K| \quad (22)$$

$$a_{iw}, x_w \in \{0,1\} \quad \forall w \in \Omega \quad (23)$$

All the equations expressed for the case of maximizing fulfillment form are the same as the mathematical model in the case of cost minimization, but they are different in the two equations. The objective function and the constraint on the maximum acceptable traveled distance are different as shown in equations (18) and (21) respectively. Although the stated mathematical model is similar to the previous one, it is used when the demand fulfillment is more important than the cost minimization.

Although route-based formulation has a simple linear structure, the calculation of feasible routes and the costs associated with each route increase significantly with increasing size. Finding routes that meet all the constraints also requires a systematic approach. For this purpose, an algorithm has been used to generate feasible routes and solve the problem accurately.

## 5. Exact algorithm

The proposed mathematical models require the calculation of feasible routes to run the model. Therefore, with the help of an algorithm, all feasible routes and their costs are calculated in MATLAB 2016. Finally, according to the objective functions, the best routes are selected using mathematical models.

The algorithm used is an approach for constructing and examining possible solutions using a specific solution structure. In this method, a solution structure is considered for the possible space of the problem. Therefore all possible solutions can be produced. This solution structure includes the customer assignment and the customer visit sequence for each vehicle. If each vehicle visits the $I$ node and $K$ represents the number of the available vehicle, then a solution string that can generate all possible routes in $I + K - 1$ length. In this solution string, the numbers $\{1,2, \ldots, I\}$ represent the customer node number, and characters * are used as separators in the solution string. For example, consider a problem with 8 customers and two vehicles. A solution string can exist as follows:

| 7 | 2 | 5 | 3 | 1 | * | 8 | 4 | 6 |

Depending on the solution string, there are two vehicles, the character * is used as the separator. Therefore, the interpretation of the solution string will be as follows:

- Route of vehicle number one: 7→2→5→3→1
- Route of vehicle number two: 8→4→6

Using the solution string, all possible routes can be generated, and in the next step, each solution must be evaluated based on the constraints. Due to the type of solution structure, several constraints, such as the number of available vehicles and the visiting of all customers, are automatically observed. Therefore, every solution should only be evaluated based on capacity constraint and the success rate constraint. Also, due to the clarity of the routes, calculating the success rate will be simple. Of course, the solving approach based on route production will not be suitable for solving large-scale problems due to the significant increase in the number of routes.

According to the explanations, the solving approach using the route generation algorithm will be as follows:

- First, generate all possible routes based on the proposed solution structure.
- Evaluate the vehicle capacity constraints and the minimum success rate (or maximum traveled distance).
- For all feasible solutions, calculate the value of the objective function (distance traveled or total success rate).
- Then select the best solution from the generated solutions based on the value of the objective function.

The proposed approach always provides the optimal solution due to the search of the entire solution space. But generating routes and evaluating various constraints requires time. Therefore, other approaches are suggested in the following sections.

The robust vehicle routing problem in the presence of a proactive attacker is evaluated in two cases. In the first case, the main goal is to minimize the cost of supply (distance traveled). Maximizing the probability of supply or demand fulfillment is the main goal of the second case. To evaluate the solution approach, the sample data in Table 1 have been used. As mentioned in the problem definition, to solve the problem, data on the success rate is needed. The success rate values for all arcs are shown in Table 6.

Table 6. the success rate for passing network arcs (the success rate of supply)

| Node number | 0 | 1 | 2 | 3 | 4 | 5 | 6 | 7 | 8 |
|---|---|---|---|---|---|---|---|---|---|
| 0 | 1 | 0.961 | 0.949 | 0.941 | 0.952 | 0.948 | 0.982 | 0.949 | 0.972 |
| 1 | 0.961 | 1 | 0.937 | 0.932 | 0.963 | 0.963 | 0.985 | 0.987 | 0.976 |
| 2 | 0.949 | 0.937 | 1 | 0.982 | 0.956 | 0.952 | 0.937 | 0.973 | 0.962 |
| 3 | 0.941 | 0.932 | 0.982 | 1 | 0.966 | 0.987 | 0.955 | 0.962 | 0.970 |
| 4 | 0.952 | 0.963 | 0.956 | 0.966 | 1 | 0.933 | 0.981 | 0.975 | 0.986 |
| 5 | 0.948 | 0.963 | 0.952 | 0.987 | 0.933 | 1 | 0.980 | 0.977 | 0.978 |
| 6 | 0.982 | 0.985 | 0.937 | 0.955 | 0.981 | 0.980 | 1 | 0.958 | 0.984 |
| 7 | 0.949 | 0.987 | 0.973 | 0.962 | 0.975 | 0.977 | 0.958 | 1 | 0.920 |
| 8 | 0.972 | 0.976 | 0.962 | 0.970 | 0.986 | 0.978 | 0.984 | 0.920 | 1 |

The implementation results of the exact algorithm in the case of cost minimization for different values of minimum success rate can be seen in Table 7. Figure 3 also shows the chart of the objective function in various sample number for different values of the success rate ($\alpha$).

Table 7. Results of the exact algorithm for the cost minimization form

| Sample number | Objective function ($Z$) | Total success rate in supply | Time (s) | Minimum expected success rate ($\alpha$) |
|---|---|---|---|---|
| 1 | 211.25 | 7.42 | 119.5 | 0 |
| 2 | 211.25 | 7.42 | 126.5 | 20 |
| 3 | 211.25 | 7.42 | 127.1 | 50 |
| 4 | 211.25 | 7.42 | 125.1 | 80 |
| 5 | 211.25 | 7.42 | 125.1 | 82 |
| 6 | 211.25 | 7.42 | 124.6 | 85 |
| 7 | 237.37 | 7.50 | 125.3 | 87 |
| 8 | 248.55 | 7.55 | 123.9 | 90 |
| 9 | 310.29 | 7.61 | 124.5 | 92 |
| 10 | 310.29 | 7.61 | 124.9 | 92.2 |
| 11 | infeasible | infeasible | 123.8 | 92.3 |

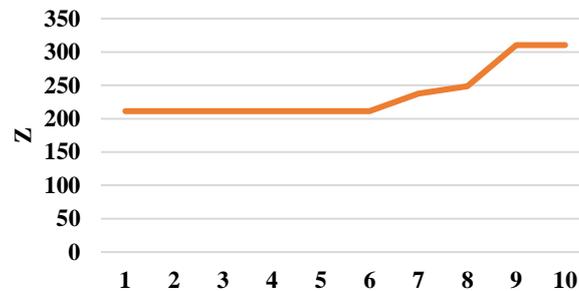

Figure 3. objective function chart

According to the objective function value obtained for the different success rates, it can be concluded that increasing the probability of expected success to 85% does not affect the objective function. The objective function value with a minimum probability of success rate of less than 87% is equal to the objective function value of classical capacitated vehicle routing problem in the optimal state (211.25). Also, as the probability of success in supply increases, the objective function value increases so that more reliable routes can be obtained. The success rate in supply is directly related to the amount of distance traversed by vehicles, and if the success rate increases, the amount of cost of vehicles will increase. Also, due to the constant process of algorithm, time-consuming for various sample data of the same size not much different from each other.

To analyze the results, the optimal solution for sample number 8 is evaluated. In this example, there are eight demand points and a central depot. The objective function value is 211.25 in non-attacking network. If the minimum supply success rate for each node is set at 90%, the objective function value will be 248.55. To increase the total success rate, vehicles have to travel longer distances. The routes are shown in Figure 4. Table 8 represents details of the optimal solution.

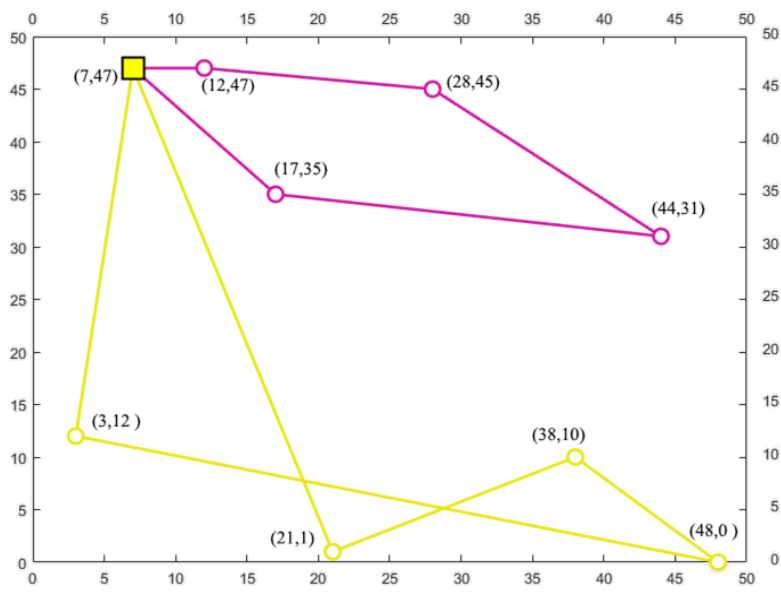
Figure 4. representation of the solution problem 8 in the case of cost minimization

Table 8. optimal solution detail for the problem 8

| Route number | Route | Used capacity | Available capacity | Objective function | $\alpha$ | Node number | Node success rate (%) |
|---|---|---|---|---|---|---|---|
| 1 | 0→6→5→1→7→0 | 157 | 180 | 248.55 | 90 | 6 | 98.2 |
| | | | | | | 5 | 96.2 |
| | | | | | | 1 | 92.7 |
| | | | | | | 7 | 91.5 |
| 2 | 0→8→4→3→2→0 | 151 | 180 | | | 8 | 97.2 |
| | | | | | | 4 | 95.8 |
| | | | | | | 3 | 92.6 |
| | | | | | | 2 | 90.9 |

When the success rate is considered, the main goal is to increase the reliability of supply and reduce costs to the next level. In this section, the proposed algorithm is used to evaluate the fulfillment maximization mode. The results of the implementation of the proposed algorithm in maximization form for different amounts of expected costs are shown in Table 9.

Table 9. Results of the exact algorithm for maximizing demand fulfillment

| Sample number | Objective function (Z) | Total cost | Time | Total expected cost ($\beta$) |
|---|---|---|---|---|
| 1 | 7.62 | 337.76 | 123.1 | 400 |
| 2 | 7.62 | 337.76 | 126.1 | 350 |
| 3 | 7.61 | 310.29 | 126.1 | 330 |
| 4 | 7.60 | 299.65 | 126.0 | 310 |
| 5 | 7.60 | 299.65 | 127.9 | 300 |
| 6 | 7.58 | 279.17 | 125.3 | 280 |
| 7 | 7.56 | 255.79 | 125.8 | 260 |
| 8 | 7.50 | 237.37 | 125.7 | 240 |
| 9 | 7.44 | 212.89 | 131.7 | 230 |
| 10 | 7.44 | 212.89 | 131.9 | 220 |
| 11 | infeasible | infeasible | 130.1 | 210 |

The results show that there is no significant difference in the solution time for different samples problem of the same size. Also, the objective function value is directly related to the expected total cost. When the maximum acceptable cost is reduced, the total success rate will be decreased. It should be noted that in sample number 11, where the maximum accepted cost is less than the optimal solution in the normal case (210 <211.25), there will be no feasible solution.

To examine the results, sample number 10 has been selected. There are eight demand points and a central depot. In this case, the maximum total acceptable cost is 220. Accordingly, the objective function value is 7.44 and the total distance traveled is 212.89. Therefore, to travel shorter distances, routes with more probability of attacking are selected. The route depicted in the figure 5 and the optimal solution is shown in Table 10.

Table 10. optimal solution detail for the problem 10

| Route number | Route | Used capacity | Available capacity | Objective function | β | Node number | Node success rate (%) |
|---|---|---|---|---|---|---|---|
| 1 | 0→6→1→7→0 | 129 | 180 | 212.89 | 220 | 6 | 98.2 |
|  |  |  |  |  |  | 1 | 96.7 |
|  |  |  |  |  |  | 7 | 95.5 |
| 2 | 0→8→2→3→4→5 | 179 | 180 |  |  | 8 | 97.2 |
|  |  |  |  |  |  | 2 | 93.5 |
|  |  |  |  |  |  | 3 | 91.8 |
|  |  |  |  |  |  | 4 | 88.7 |
|  |  |  |  |  |  | 5 | 82.8 |

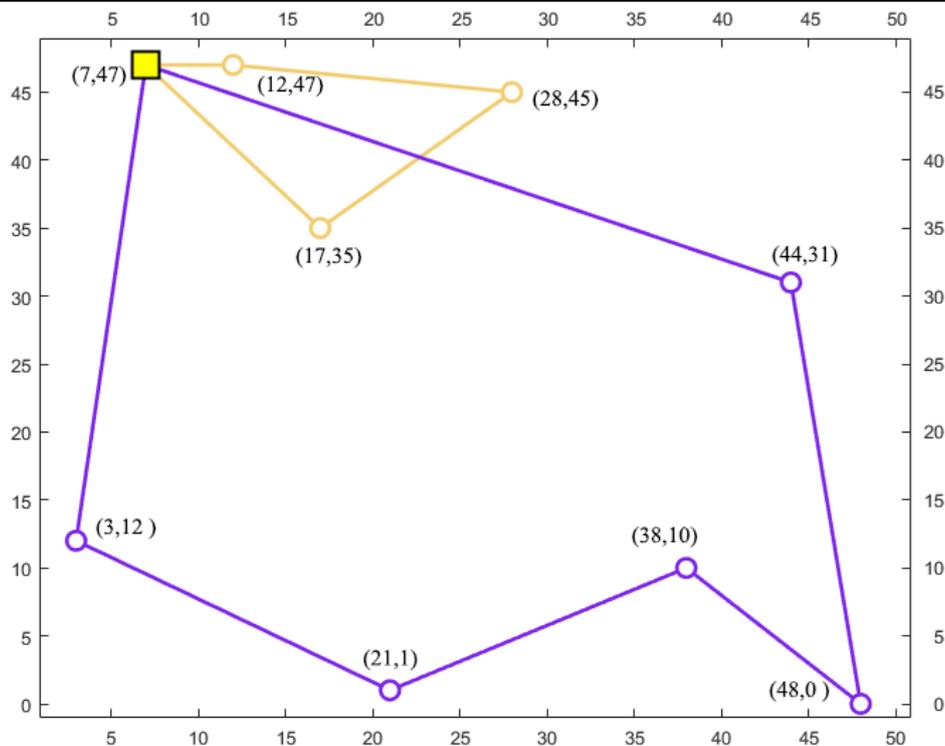

Figure 5 representation of the solution problem 10 in the case of fulfillment minimization.

The vehicle routing problem is NP-hard [15], and finding the exact solution will be time-consuming and complex. The proposed algorithm will have a significant implementation time due to the route-based structure. These conditions will be more difficult for large-scale problems. Although the exact algorithm can solve the problem accurately, using this approach to solve large-scale problems is not possible. Therefore, it is better to use faster solution methods such as meta-heuristic algorithms. Meta-heuristic algorithms are random-based algorithms that are used to find optimal or near-optimal solutions. The structure of all meta-heuristic algorithms must have two basic features. The first feature is the use of appropriate local search methods, and the second is the ability to escape local-optimum. In this research, to increase the speed of solving large-scale problems, a simulated annealing algorithm (SA) has been used.

## 6. Simulated Annealing Algorithm

An optimization problem is a problem of finding the best solution among all possible solutions. NP-hard problems are a group of optimization problems in which it is unlikely to exist a polynomial algorithm to solve them. So, it is difficult and very timely costing to manage and find the optimal solution using exact methods due to the increasing number of feasible solutions. One of the practical algorithms is the simulated annealing algorithm. Simulated annealing algorithms are often used to estimate global optimization for large-scale problems. Also, this algorithm is usually used when the search space is discrete [16].

In materials science, annealing is the thermal processing during which the physical and sometimes the chemical properties of a material change. This is done to increase the malleability and reduce the hardness of the material. During this process, the metal is first heated, then kept at a certain temperature, and finally, slowly cooled. As the metal heats up, the molecules move freely in any direction. This freedom decreases as the material gradually cools. If the cooling process is slow enough, it can be ensured that the heat energy is evenly distributed throughout the body. Therefore it has the best crystal structure that is symmetrical and resistant. The solution heats up well; then, the variation range slowly decreases to find the optimal or near-optimal solution.

Due to the applications of the simulated annealing algorithm in VRP, this approach has been used to solve the robust vehicle routing problem in the presence of a proactive attacker. Although in the previous sections, two different forms of the problem (minimizing costs and maximizing the total success rate) were evaluated, the SA algorithm can solve both types of problems.

The general structure of the simulated annealing algorithm has several main parts that act as a general structure. These parts are interconnected and the capability and efficiency of each section improve the entire structure. The SA algorithm has three main parts which are related to the process of starting, searching for the solution space, selecting the good solutions, and finally accepting the best-searched solution. These parts are defined as follows:

- Appropriate initial solution
- Initial temperature and a temperature reduction structure to avoid local optimum
- Local search structures to find the good solutions

Using appropriate methods to increase the efficiency of the parts, increases the structure efficiency. How to define these items is discussed below.

## 6.1. Initial temperature, temperature reduction structure and initial solution

One of the main parts of the SA algorithm is to calculate the appropriate initial temperature and the temperature reduction approach in each iteration. To use the SA, it is necessary to have an appropriate initial solution and then apply a local search structure to find new solutions. To find better solutions the search space is transferred to another search space with a certain probability. If the objective function value improves, the solution is accepted. Otherwise, this solution will be accepted with a certain probability. The probability of accepting worse solutions will be calculated based on two factors: the amount of solution difference and the current temperature. This approach leads to the acceptance of different solutions and the escape from the local optimum. The probability of accepting the worse solutions decreases with a temperature reduction. Therefore, the overall performance of the algorithm depends entirely on the initial solution, initial temperature, temperature reduction structure, and how the neighborhood search is defined [17].

In this research, to define the initial temperature, a random solution structure has been used. First, several answers are generated randomly and then the initial temperature is determined based on the objective function values. How to calculate the initial temperature ($T_0$) represented based on the equations 24-26.

$$f_{max} = max_{i=1}^{n}\{f_i\} \qquad (24)$$
$$f_{min} = min_{i=1}^{n}\{f_i\} \qquad (25)$$
$$T_0 = f_{min} + 0.1 * (f_{max} - f_{min}) \qquad (26)$$

Equation 24 is defined to calculate the maximum objective function value. Equation 25 shows the minimum objective function value among the randomly generated solutions. Finally, the initial temperature is calculated through Equation 26.

Note that the initial solution is selected from randomly generated solutions to calculate the initial temperature, in other words, the initial solution is equal to the best generated solution among the random solutions.

$$S_{initial} = \{X | x \in X_{f_{min}}\} \qquad (27)$$

The temperature reduction structure in the SA algorithm should be designed in such a way that at the beginning of the algorithm the probability of accepting worse solutions is high. After a proper search of the solution space, the probability of acceptance is reduced. For this purpose, a linear temperature reduction structure is used. The temperature reduction equation is shown below:

$$T_i = \delta . T_{i-1} \qquad (28)$$

In equation 28 the new temperature value of each iteration is a coefficient of the previous temperature. Noted that to reduce the temperature structure, a number between zero and one is selected ($0 < \delta < 1$).

## 6.2. Local search structure

In meta-heuristic algorithms, local search structures are used to search in a specific neighborhood and find the optimal local solution. These operators must be designed based on the structure to achieve better solutions. In this research, three operators include shift, exchange and 2-opt have been used. These operators are designed according to the structure of the solution space and increase the possibility of finding different and better solutions [18].

- Shift operator: a demand point in the route is randomly selected each time. Then it is randomly placed among other possible points (other possible positions). Figure 6 shows an overview of the operator.

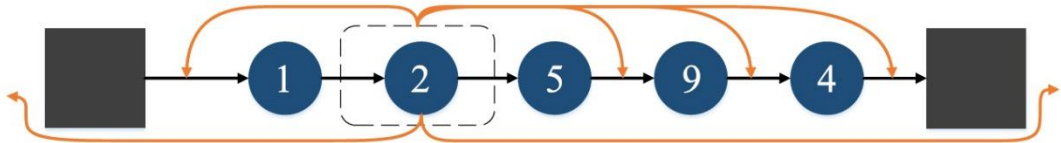

Figure 6. shift operator structure

- Exchange operator: The exchange operator is used to change the position of two demand points in the customer visit sequence. In fact, at each stage, two demand points are randomly selected and their positions are swapped. Figure 7 shows an overview of the exchange operator.

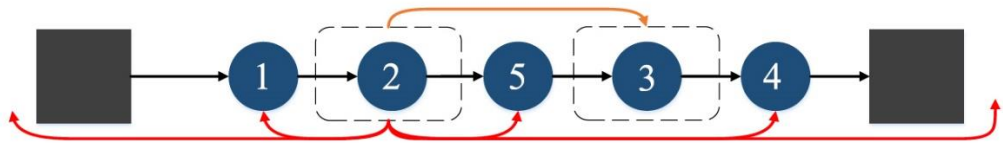

Figure 7. exchange operator structure

- 2-opt operator: two arcs are randomly removed from a route and two new arcs replace the previous ones. An overview of the 2-opt operator is shown in Figure 8.

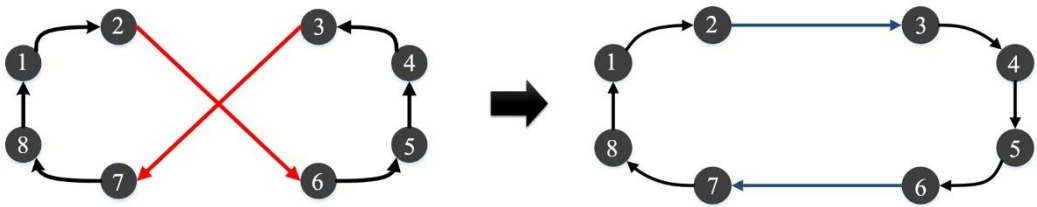

Figure 8. 2-opt operator structure

After defining the details, the main structure of the simulated annealing algorithm can be defined. After calculating the appropriate initial solution by generating random solutions, the initial solution enters the algorithm process. This process has two main loops. The initial loop specifies the total number of iterations. The second loop is the local search structure. In the local search loop, better solutions will be accepted, and worse solutions will be accepted or rejected based on a specific probability. After completing the inner loop operation, the temperature reduction structure is used and the temperature is updated. This process is repeated until the end of the main loop iterations. The main structure of the simulated annealing algorithm is represented in the pseudo-code of Figure 9.

| | |
|---|---|
| *Simulated Annealing Structure* | |

Select an initial solution $S$
Select the temperature cooling schedule $T_k = \alpha . T_{k-1}$
Select an initial temperature $T = T_0$
Select a repetition schedule, $L$, that defines the number of main loop iterations
Select a repetition schedule, $K$, that defines the number of iterations executed at each temperature, $T_k$
**Repeat** until $l = L$
Set repetition counter $l = 0$
    **Repeat** until $k = K$
    Set repetition counter $k = 0$
    Generate a solution $S'$ with **local search**
    **Calculate** $\Delta_{S,S'} = f(S') - f(S)$
    **If** $\Delta_{S,S'} \leq 0$, **then** $S \leftarrow S'$
    **If** $\Delta_{S,S'} \leq 0$, **then** $S \leftarrow S'$ with probability $exp(-\Delta_{S,S'}/T_k)$
    $k \leftarrow k+1$
$T_k = \alpha . T_{k-1}$
$l \leftarrow l+1$

Figure 9. proposed simulated annealing algorithm pseudo-code

## 7. Numerical Results

Evaluating the efficiency of the simulated annealing algorithm requires various examples. For this reason, the meta-heuristic algorithm is first applied to the previous sample data and the results are compared with the exact approach. Then some large-scale problems are generated and evaluated using the proposed algorithm.

In the previous sections, the problem was considered in cost minimization and fulfillment maximization forms. But the proposed SA algorithm is defined in the minimization mode. Therefore, to solve the second case, i.e. maximizing the total success rate, it is necessary to make small changes in the calculation of the objective function. Therefore, when calculating the objective function of the problem in maximization mode, the value of the objective function is calculated negatively.

$$f(X) = Max(Z) = Min(-Z) \qquad (29)$$

Also, due to the random generation of initial solutions, a feasible solution may not be generated at first. Therefore, to achieve the feasible solutions, for not observing any of the constraints, a penalty function is considered which is added to the objective function. This moves new solutions to the feasible space. In this case, if the constraint is $AX \leq b$, the penalty function, and the general objective function are calculated as follows.

$$f_p(X) = M . \max((AX - b), 0) \qquad (30)$$
$$f_t(X) = f(X) + f_p(X) \qquad (31)$$

The proposed algorithm was evaluated on 22 sample data. 11 samples are for cost minimization form and 11 samples are for maximizing total success rate, and the computational results are shown in Table 11.

Table 11. Results of simulated annealing algorithm compared to the exact algorithm

| Sample number | Optimal objective function (Z*) | Time$_e$ (s) | Z$_{SA}$ | Time$_{SA}$ (s) | Gap (%) |
|---|---|---|---|---|---|
| 1 | 211.25 | 119.5 | 211.25 | 5.9 | 0 |
| 2 | 211.25 | 126.5 | 211.25 | 6.2 | 0 |
| 3 | 211.25 | 127.1 | 211.25 | 5.9 | 0 |
| 4 | 211.25 | 125.1 | 211.25 | 6.1 | 0 |
| 5 | 211.25 | 125.1 | 211.25 | 6.2 | 0 |
| 6 | 211.25 | 124.6 | 211.25 | 6.0 | 0 |
| 7 | 237.37 | 125.3 | 237.37 | 6.3 | 0 |
| 8 | 248.55 | 123.9 | 248.55 | 6.3 | 0 |
| 9 | 310.29 | 124.5 | 310.29 | 6.2 | 0 |
| 10 | 310.29 | 124.9 | 310.29 | 6.3 | 0 |
| 11 | Infeasible | 123.8 | Infeasible | 6.4 | 0 |
| Avg (1-11) | 237.4 | 124.6 | 237.4 | 6.2 | 0 |
| 12 | 7.62 | 123.1 | 7.62 | 6.2 | 0 |
| 13 | 7.62 | 126.1 | 7.62 | 6.3 | 0 |
| 14 | 7.61 | 126.1 | 7.61 | 6.4 | 0 |
| 15 | 7.60 | 126.0 | 7.60 | 6.2 | 0 |
| 16 | 7.60 | 127.9 | 7.60 | 6.1 | 0 |
| 17 | 7.58 | 125.3 | 7.58 | 6.0 | 0 |
| 18 | 7.56 | 125.8 | 7.56 | 6.5 | 0 |
| 19 | 7.50 | 125.7 | 7.50 | 6.5 | 0 |
| 20 | 7.44 | 131.7 | 7.44 | 6.6 | 0 |
| 21 | 7.44 | 131.9 | 7.44 | 6.8 | 0 |
| 22 | Infeasible | 130.1 | Infeasible | 6.8 | 0 |
| Avg (12-22) | 7.56 | 127.2 | 7.56 | 6.4 | 0 |

According to Table 11 in all the sample data, the proposed meta-heuristic algorithm has achieved the optimal solution. Moreover, the proposed meta-heuristic algorithm is much faster than the exact algorithm. Therefore, good performance can be expected from the proposed algorithm. Note that the gap is calculated using the equation $((Z_{SA} - Z^*)/Z^*) * 100$.

The efficiency of the proposed meta-heuristic algorithm on a small scale was compared with the optimal solution and showed the speed and capability. To evaluate the proposed algorithm on various scales, 6 examples are generated and the proposed algorithm is evaluated. The method of producing numerical examples is shown in Table 12.

Table 12. data generation structure

| Symbol | Description | Production Function |
|---|---|---|
| $(x, y)$ | Coordinates of demand points | $U(0,100)$ |
| $d_i$ | Demand | $U(20,50)$ |
| $p_{ij}$ | The success rate for crossing arc $(i, j)$ | $U(0.9, 0.99)$ |
| $Q_k$ | Vehicle capacity | $\left\lceil 1.02 \frac{\sum_i d_i}{|K|} \right\rceil$ |

The results of the proposed SA algorithm are shown on 6 sample problems in Table 13. 6 numerical examples from 10 to 100 demand points and from 2 to 5 vehicles have been evaluated. The total cost (traveled distance) and the total supply success rate are denoted by C and φ, respectively. To more accurately evaluate, the total cost and the total success rate with a non-interdicted network ($\alpha = 0$) have been implemented. In this case, the problem becomes a classical CVRP.

Table 23. Results of simulated annealing algorithm un different scale problem

| #N | N | |K| | $C_{\alpha=0}$ | $\varphi_{\alpha=0}$ | T(s) | α | $C_\alpha$ | $\varphi_\alpha$ | T(s) | β | $C_\beta$ | $\varphi_\beta$ | T(s) |
|---|---|---|---|---|---|---|---|---|---|---|---|---|---|
| 1 | 10 | 2 | 453.78 | 8.438 | 72.1 | 75% | 453.78 | 8.6048 | 73.6 | 500 | 485.46 | 9.0193 | 38.1 |
| 2 | 20 | 2 | 439.43 | 15.1083 | 74.2 | 75% | 627.22 | 17.3564 | 111.9 | 700 | 697.19 | 17.8699 | 59.5 |
| 3 | 30 | 3 | 566.94 | 21.5225 | 120.8 | 75% | 927.07 | 25.9148 | 136.8 | 1000 | 998.02 | 26.4863 | 81.9 |
| 4 | 40 | 3 | 709.22 | 28.2805 | 153.7 | 75% | 1333.03 | 33.955 | 162.1 | 1400 | 1383.99 | 35.3427 | 94.8 |
| 5 | 50 | 4 | 832.84 | 35.7988 | 160.9 | 75% | 1763.89 | 42.5133 | 162.3 | 1900 | 1885.49 | 44.2569 | 97.6 |
| 6 | 100 | 5 | 1045.23 | 63.183 | 178.4 | 75% | 2760 | 86.3254 | 198.1 | 2800 | 2777.57 | 90.1156 | 118.2 |

Chart of total cost and the total success rate in the stable network (classical CVRP), minimization of the total cost in the unstable network, and maximization of demand fulfillment in the unstable network are shown in Figures 10 and 11, respectively.

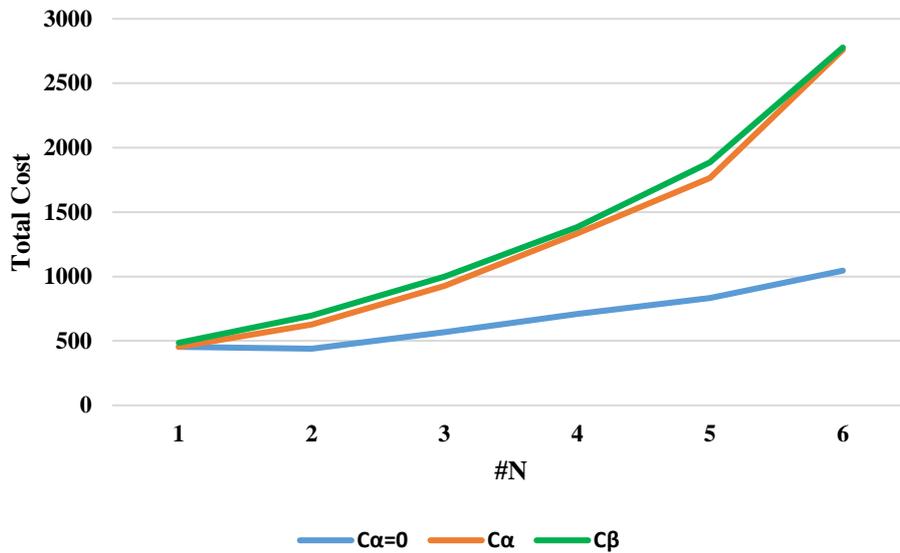

Figure 10. total cost for the stable network, cost minimization, and fulfillment maximization form

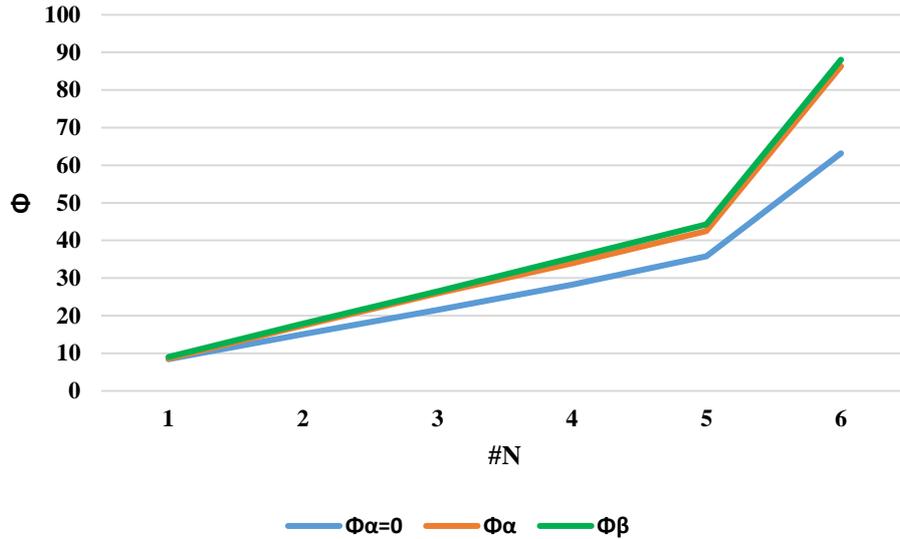

Figure 11. total success rate in the stable network, cost minimization, and fulfillment maximization form

According to the results, instability in the network can increase the cost. In other words, longer routes are usually chosen to increase the total success rate. Therefore, routes with more reliability are generally longer distances. The results also show the efficiency of the proposed algorithm in solving problems of different sizes.

## 8. Conclusions

The application of the vehicle routing problem in the field of the supply chain is significant. Also, due to the importance of instability in the network routing problem, a hybrid problem is defined as the robust vehicle routing problem in the presence of a proactive attacker. Due to uncertain conditions when vehicles pass through the arcs, there is a possibility of attacking in the route, vehicle breakdown, and incapability to fulfill demand. The possibility of attacking the routes in war conditions or the instability of the routes during natural disasters are the causes of instability. The concept of instability refers to the conflict between attackers and defenders in a network. In this study, the probability of attacking is estimated by the network user and then the vehicle routing problem on the unstable network is defined.

In this study, the classical capacitated vehicle routing problem was first evaluated. In this case, to calculate the probability of attacking, the sequence of customer visits in each route is required. For this reason, an extensive model of the vehicle routing problem was used to make it easier to calculate the probability of attacking. Then, with the help of the extended model, two modes of cost minimization and maximization of demand fulfillment were defined. In the case of cost minimization, the distance traveled by vehicles is minimized when there is a minimum success rate for each node. In the maximization form, increasing the total success rate is considered the main goal.

For evaluation, an exact algorithm based on the production of feasible routes was presented. In this algorithm, only routes that conform to the constraints are generated and investigated. Then, by examining a numerical example, the effects of the parameters on the results were determined. When there is an unstable network, the total cost to satisfy customers will increase. Increasing the reliability of crossing routes will increase the total cost. Also, due to the incremental computational cost of using the exact algorithm, a simulated annealing algorithm was proposed. In this algorithm, by using the appropriate initial temperature, temperature reduction structure, and local search operators, an extensive search is created in the solution space, which leads to finding optimal or near-optimal solutions. To evaluate the proposed meta-heuristic algorithm, six numerical problems on the different scales were used. Then, three modes of the problem were evaluated for stable conditions, cost minimization with instability, and fulfillment maximization with instability using the proposed algorithm. The results show the efficiency of the algorithm in different sizes. Also, comparing the results of the exact algorithm and the proposed meta-heuristic algorithm shows the appropriate capability of the method to finding optimal solutions.